\documentclass{amsart}
\usepackage[alphabetic]{amsrefs}
\usepackage{amsmath, amsthm, graphicx, amssymb,verbatim,pst-all,color}
\usepackage[all]{xy}
\ifx\JPicScale\undefined\def\JPicScale{1.0}\fi
\unitlength \JPicScale mm
{\theoremstyle{plain}
   
   \newtheorem{theorem}{Theorem}[section]
    \newtheorem{theorem*}{Theorem*}[section]
   \newtheorem{proposition}[theorem]{Proposition}     
   \newtheorem{lemma}[theorem]{Lemma}
   
   \newtheorem{corollary}[theorem]{Corollary}

}
{\theoremstyle{definition}

   \newtheorem{definition}[theorem]{Definition}
   \newtheorem{remark}[theorem]{Remark}

}
\newcommand{\quotient}{/\hspace{-1.2mm}/}

\newcommand{\M}{\overline{\mathcal{M}}}
\newcommand{\SL}{\operatorname{SL}}

\newcommand{\Amp}{\operatorname{Amp}}
\newcommand{\Nef}{\operatorname{Nef}}

\newcommand{\dist}{\operatorname{dist}}
\newcommand{\HH}{\operatorname{H}}
\newcommand{\N}{\operatorname{N}}
\newcommand{\MM}{\operatorname{M}}
\newcommand{\Proj}{\operatorname{Proj}}

\newcommand{\mc}[1]{\mathcal{#1}}
\newcommand{\ov}[1]{\overline{#1}}
\newcommand{\op}[1]{\operatorname{#1}}
\newcommand{\ovop}[1]{\ov{\op{#1}}}

\newcommand{\ovmc}[1]{\ov{\mc{#1}}}
\newcommand{\sL}{\mathfrak{sl}}

\numberwithin{theorem}{section}

\begin{document}
\title[The cone of type A, level one conformal blocks divisors]{the cone of type A, level one \\ conformal blocks divisors}
\author{Noah Giansiracusa and Angela Gibney}
\maketitle


\begin{abstract}
We prove that the type $A$, level one, conformal blocks divisors on $\ovop{M}_{0,n}$ span a finitely generated, full-dimensional subcone of the nef cone.  Each such divisor induces a morphism from $\ovop{M}_{0,n}$, and we identify its image as a GIT quotient parameterizing configurations of points supported on a flat limit of Veronese curves.  We show how scaling GIT linearizations gives geometric meaning to certain identities among conformal blocks divisor classes.   This also gives modular interpretations, in the form of GIT constructions, to the images of the hyperelliptic and cyclic trigonal loci in $\ovop{M}_g$ under an extended Torelli map.
\end{abstract}

\section{Introduction}
The spaces $\M_{g,n}$ occupy a special position in the landscape of algebraic geometry.  As moduli spaces, they give insight into smooth curves and their degenerations.  As varieties, they are fertile ground for the application of rapidly developing techniques in birational geometry and the minimal model program \cite{HuKeel,GKM,HassettWeighted, simpsonlogcanonical,HHFlip,HassettHyeon,CastravetTevelev,HyeonLee,BCHM}.  We study the role played by certain vector bundles on $\ovop{M}_{0,n}$ arising in conformal field theory, especially with regard to the Mori dream space (MDS) conjecture of Hu and Keel.

The MDS conjecture states that $\ovop{M}_{0,n}$ is a Mori dream space, and hence that the cone of nef divisors is finitely generated and every nef divisor is semi-ample \cite{HuKeel}.  There are nef divisors on $\ovop{M}_{g,n}$ that are not semi-ample when $g>0$ \cite{KeelAnnals2}, but the MDS conjecture would imply that for any $g\ge 0$, there are only finitely many morphisms with connected fibers from $\ovop{M}_{g,n}$ to a projective variety \cite{HuKeel,GKM}.  This motivates the problem of classifying all such morphisms and describing their images.  We carry out this program for the set 
$\mathcal{CB}(A,1)$ of level one conformal blocks divisors on $\ovop{M}_{0,n}$ with type $A$ root system by (1) identifying the images of all induced morphisms, (2) showing that the cone it spans is finitely generated, and (3) establishing geometrically meaningful identities among the cone's generators.
 
\subsection{Finite generation}
Fakhruddin has very recently unlocked a treasure chest containing a potentially infinite collection of semi-ample divisors on $\ovop{M}_{0,n}$ \cite{Fak}.  Each of these {\textit{conformal blocks divisors}},  $\mathbb{D}^{\mathfrak{g}}_{\ell, w}$, is the first Chern class of a vector bundle 
$\mathbb{V}^{\mathfrak{g}}_{\ell, w}$. These bundles, defined on the stacks $\ovmc{M}_{g,n}$ for each $g$ and $n$, are cooked up from three ingredients: a simple lie algebra $\mathfrak{g}$, a positive integer $\ell$, and an $n$-tuple of dominant integral weights $w$ for $\mathfrak{g}$ of level $\le \ell$ \cite{TUY}.  In the '90s a great deal of work went into proving the Verlinde formula and showing that there are canonical identifications of the fiber of $\mathbb{V}^{\mathfrak{g}}_{\ell, w}$ over a point $X= (C,p_1,\ldots,p_n) \in \mc{M}_{g,n}$ with a space of generalized theta functions on $C$ \cite{BertramSzenes, BeauvilleLaszlo, Faltings, KNR, Pauly, LaszloSorger}.  For example, $\mathbb{V}^{\sL_r}_{ \ell, w}|_X \cong \Gamma(\operatorname{SU}_C(r, w), \mathcal{L})$, where $\operatorname{SU}_C(r, w)$ is the moduli stack of rank $r$ vector bundles on $C$ with trivial determinant and parabolic structure at the $p_i$ given by $w$, and $\mathcal{L}$ is a line bundle determined by $\ell$ and $w$ \cite{Pauly}.  

Fakhruddin has refocused attention on these bundles by showing that for $g=0$, the $\mathbb{V}^{\mathfrak{g}}_{\ell, \omega}$ are globally generated, so their first Chern classes $\mathbb{D}^{\mathfrak{g}}_{\ell, \omega}$ are semi-ample divisors \cite[Lemma 2.5]{Fak}.   This both lends support to, and casts a shadow over, the MDS conjecture.  On the one hand, it provides a seemingly infinite collection of nef divisors on $\ovop{M}_{0,n}$ which, as predicted by the conjecture, are all semi-ample.  On the other hand, a number of these divisors span extremal rays of the nef cone \cite{ags, agss}, and the conjecture says that the nef cone should have only finitely many extremal rays.  We prove that, for a large class of conformal blocks bundles, this potential infinitude of extremal rays cannot occur:

\begin{theorem}\label{thm:finite}
The cone in $\N^1(\ovop{M}_{0,n})$ spanned by the divisor classes $$\mathcal{CB}(A,1) := \{\mathbb{D}^{\sL_m}_{1, w} : w \mbox{ is arbitrary and }m \ge 2\}$$ is finitely generated.  In particular, these divisors yield only finitely many extremal rays of $\Nef(\ovop{M}_{0,n}).$
\end{theorem}  

This cone is full-dimensional, since the elements of $\mathcal{CB}(\sL_2,1)$ form a basis for $\op{Pic}(\ovop{M}_{0,n})$ \cite[Theorem 4.3]{Fak}.  Our approach to Theorem \ref{thm:finite} is to show first that this cone is spanned by a finite number of cones arising in a GIT construction.  These GIT cones, generalizing the one introduced in \cite{AS08}, are then shown to be finitely generated by a direct application of variational results of Thaddeus and Dolgachev-Hu \cite{Tha96,DH98}. 

\subsection{Induced morphisms}
Each element of $\mathcal{CB}(A,1)$ gives rise to a morphism from $\ovop{M}_{0,n}$ to a projective variety.  We show that these morphisms factor as the composition of a forgetful map with a birational morphism to a GIT quotient parameterizing configurations of points that lie on a flat limit of Veronese curves---what is called in \cite{Giansiracusa} a quasi-Veronese curve.  This birational map further factors through a Hassett moduli space $\ovop{M}_{0,\vec{c}}$ of weighted pointed curves.  More precisely, if $V_{d,k}$ denotes the locus of configurations of $k$ points in $\mathbb{P}^d$ lying on a quasi-Veronese curve, then $\SL_{d+1}$ naturally acts on $V_{d,k}$, and we have the following:

\begin{theorem}\label{thm:morphism} The morphism given by a multiple of the semi-ample divisor $\mathbb{D}^{\sL_m}_{1,w}$ is equal to the composition \begin{equation}\label{eq:composition}\ovop{M}_{0,n} \overset{\pi}{\longrightarrow} \ovop{M}_{0,k} \longrightarrow \ovop{M}_{0,\vec{c}} \longrightarrow V_{d,k}\quotient_{\vec{c}}\SL_{d+1},\end{equation} where $\pi$ drops the points with zero weight, and $\vec{c}$ corresponds to the nonzero weights. \end{theorem}  

The relation between $w,m$, and $d$ is explained in \S\ref{section:RelGIT}.  This map to $V_{d,k}\quotient_{\vec{c}}\SL_{d+1}$ was introduced in \cite[Theorem 1.1]{Giansiracusa}, and it was shown that the above result holds for elements of $\mathcal{CB}(\mathfrak{sl}_n,1)$ with $S_n$-invariant weights \cite[Theorem 1.2]{Giansiracusa}.  Theorem \ref{thm:morphism} vastly generalizes this result by showing that the structural decomposition of these induced morphisms does not depend on the Lie algebra or weights, it only depends on the level and the type of the root system.

Our technique for proving Theorem \ref{thm:morphism} is to compare intersection numbers with 1-strata of the boundary (called $\op{F}$-curves, cf. \cite[Theorem 2.2]{GKM}), since they span the 1-cycles on $\ovop{M}_{0,n}$.  This leads directly to a converse result: every morphism as in (\ref{eq:composition}) is given by an element of $\mathcal{CB}(A,1)$.  This builds a bridge between conformal blocks divisors and GIT quotients that, as we will see, allows information from one side to travel to the other, yielding new results in both realms.

\subsection{Divisor class identities}
The cone spanned by $\mathcal{CB}(A,1)$ is finitely generated, yet there are infinitely many type $A$ Lie algebras and possible weight vectors, so one expects many relations between the cone's generators.  An example is the following:
$$\mathbb{D}^{\mathfrak{sl}_n}_{1, (\omega_k,\ldots,\omega_k)} = \mathbb{D}^{\mathfrak{sl}_n}_{1,(\omega_{n-k},\ldots,\omega_{n-k})}.$$ This is implied by the symmetry of the Dynkin diagram for $\mathfrak{sl}_n$, and it can be verified by using Fakhruddin's formula for intersecting elements of $\mathcal{CB}(A,1)$ with $\op{F}$-curves \cite[Proposition 5.2]{Fak}.   This identity implies that the images of the induced morphisms have isomorphic normalizations: $$(V_{k-1,n}\quotient\SL_{k})^\nu\cong (V_{n-k-1,n}\quotient\SL_{n-k})^\nu.$$  This isomorphism is a manifestation of the Gale transform, so there is a geometric reason underlying the above equality of divisors classes \cite[\S6.2]{Giansiracusa}.  
 
We exhibit another divisor class identity for the cone's generators, relating elements of $\mathcal{CB}(A,1)$ with different weights \emph{and} different Lie algebras:
 
\begin{proposition}\label{introidentities}
For any  integer $k\ge 1$, one has 
$$\mathbb{D}^{\sL_m}_{1,(\omega_{c_1},\ldots, \omega_{c_n})}=\frac{1}{k}\mathbb{D}^{\sL_{km}}_{1,(\omega_{kc_1},\ldots, \omega_{kc_n})}.$$
\end{proposition}

We prove this combinatorially, using F-curve intersections, and geometrically, by showing that it corresponds to scaling GIT polarizations, via Theorem \ref{thm:morphism}.  

As a first application of Proposition \ref{introidentities}, we show that any element of $\mathcal{CB}(A,1)$ with $\op{S}_n$ invariant weights spans an extremal ray of the $S_n$-invariant nef cone.   By \cite[Theorem B]{agss}, the divisors $\{\mathbb{D}^{\sL_n}_{1, \{j,\ldots,j\}}: 2 \le j \le \lfloor \frac{n}{2} \rfloor\}$ span distinct extremal rays of $\op{Nef}(\ovop{M}_{0,n}/\op{S}_n)$, and in Corollary \ref{extremal} we use Proposition \ref{introidentities} to show that the divisors $\{\mathbb{D}^{\sL_m}_{1, \{j,\ldots,j\}}: \mbox{all } m, \ 1 \le j \le \frac{m}{2} \rfloor\}$ are proportional to these, so they also span extremal rays.  

The two types of divisor class identities described above have geometric origins, yet they are manifest as relations among generators for the conformal blocks cone.  We show that there are also situations in which conformal blocks divisor identities imply previously unknown geometric results.  For example, in Proposition \ref{SatakeImage}, using only results about divisor classes, we give GIT constructions of, and modular interpretations for, the hyperelliptic and cyclic trigonal loci in Satake's compactification of the moduli space of principally polarized abelian varieties.

\bigskip

{\em{Organization of the paper:}} In \S\ref{morphisms} we recall the morphism $\ovop{M}_{0,n} \rightarrow V_{d,n}\quotient_{\vec{x}}\SL_{d+1}$ from \cite{Giansiracusa} and prove our main technical tool, Theorem \ref{thm:ASgen}, which gives a formula for the degree on $\op{F}$-curves. This is used in \S\ref{CBGIT} to prove Theorem \ref{thm:morphism}, that elements of $\mathcal{CB}(A,1)$ induce these morphisms.  In \S\ref{finitegeneration} we prove finiteness of the cone spanned by $\mathcal{CB}(A,1)$, namely Theorem \ref{thm:finite}. In \S\ref{section:Identities} we prove the divisor class identity of Proposition \ref{introidentities} and then derive Corollary \ref{extremal} and Proposition \ref{SatakeImage}.

\bigskip

{\em{Acknowledgements:}} We thank Dan Abramovich, Najmuddin Fakhruddin, and Michael Thaddeus for their generous ideas and assistance, and we thank Maksym Fedorchuk for providing helpful feedback on an early draft.

\section{Moduli of points on quasi-Veronese curves}\label{morphisms}
The main result of this section is Theorem \ref{thm:ASgen}, which gives a simple formula for the intersection of $\op{F}$-curves with the divisors that induce the morphisms \begin{equation}\label{eq:qV}\ovop{M}_{0,n}  \overset{\varphi_{d,\vec{x}}}{\longrightarrow} V_{d,n}\quotient_{\vec{x}}\SL_{d+1},\end{equation} introduced by the first author in \cite{Giansiracusa}.  Recall that an $\op{F}$-curve on $\ovop{M}_{0,n}$ is any curve numerically equivalent to a component of the closed locus of 
points $(C,p_1,\ldots,p_n)\in \ovop{M}_{0,n}$ such that $C$ has at least $n-4$ nodes.  These curves span $\N_1(\ovop{M}_{0,n})$ by \cite{KeelTransactions} and are pictured in \cite[Theorem 2.2]{GKM}.   Each is given by a partition $\{1,\ldots,n\}=\sqcup_{i=1}^4 n_i$ with $|n_i|\ge 1$ and is denoted by $F_{n_1,\ldots,n_4}$.

\subsection{Setup}
For $1\le d \le n-3$, we denote by $U_{d,n}\subseteq(\mathbb{P}^d)^n$ the set of configurations of $n$ distinct points supported on a rational normal curve of degree $d$, so that $U_{d,n}/\text{SL}_{d+1} \cong \MM_{0,n}$.  The closure $V_{d,n} := \overline{U}_{d,n}\subseteq(\mathbb{P}^d)^n$ is the set of configurations of possibly coincident points supported on a quasi-Veronese curve, i.e., a flat limit of Veronese curves \cite[Lemma 2.3]{Giansiracusa}.  The quotients of $V_{d,n}$ by $\SL_{d+1}$, therefore, parameterize such configurations up to projectivity.

By \cite[Theorem 1.1]{Giansiracusa}, for any linearization \[\vec{x}=(x_1,\ldots,x_n)\in\Amp((\mathbb{P}^d)^n)_{\mathbb{Q}}=\mathbb{Q}^n_{>0}\] with nonempty stable locus, there is a birational morphism $\varphi_{d,\vec{x}}$ as in (\ref{eq:qV}).  The image of each $(C,p_1,\ldots,p_n)\in\ovop{M}_{0,n}$ under this morphism is obtained by first mapping $C$ to a degree $d$ quasi-Veronese curve in $\mathbb{P}^d$, and then taking the $\SL_{d+1}$-orbit of the resulting configuration of marked points.  If $C$ is smooth, then the map $C \rightarrow \mathbb{P}^d$ is simply the $d^{\text{th}}$ Veronese map, so $\varphi_{d,\vec{x}}$ induces an isomorphism on the interior: \[\ovop{M}_{0,n} \supseteq \MM_{0,n}~\widetilde{\rightarrow}~U_{d,n}/\SL_{d+1}\subseteq V_{d,n}\quotient_{d,\vec{x}}\SL_{d+1}.\]  If $C$ is nodal, then the degree $d$ map $C \rightarrow \mathbb{P}^d$ is determined by the GIT stability of the resulting configuration of points, and hence by the linearization $\vec{x}$.

\subsection{GIT polarizations}
The closure of the space of effective linearizations, restricted from $(\mathbb{P}^d)^n$, is identified with a hypersimplex \cite[Example 3.3.24]{DH98}: \[\Delta(d+1,n) := \{(x_1,\ldots,x_n)\in [0,1]^n~|~\sum_{i=1}^n x_i = d+1\}.\]  Every GIT quotient is naturally endowed with a fractional polarization.  We denote by $\mathcal{O}_{d,\vec{x}}(1)$, or simply $\mathcal{O}(1)$, the natural polarization on $V_{d,n}\quotient_{\vec{x}}\SL_{d+1}$.

The line bundle $\varphi_{d,\vec{x}}^*\mathcal{O}(1)$ on $\ovop{M}_{0,n}$ is nef, and one can ask where in the nef cone it lies as a function of $\vec{x}\in\Delta(d+1,n)$.  Since F-curve classes generate the Chow group of 1-cycles on $\ovop{M}_{0,n}$, to determine the class of a divisor it is enough to determine its degree on $\op{F}$-curves.  Alexeev and Swinarski computed this for the GIT polarizations pulled back from $(\mathbb{P}^1)^n\quotient\SL_2$ \cite[Lemma 2.2]{AS08}.  This corresponds to $d=1$ in our notation, and their formula generalizes to arbitrary $d$ as follows: 

\begin{theorem}\label{thm:ASgen}
Let $\{1,\ldots,n\}=\sqcup_{i=1}^4 n_i$ be a partition such that $|n_i| \ge 1$, and let $F_{n_1,\ldots,n_4}$ be the associated F-curve class.  For $1\le d \le n-3$ and $\vec{x}\in\Delta(d+1,n)$, set $x_{n_i} := \sum_{j\in n_i}x_j$.  Then 
\begin{equation}\label{eq:ASgen} \deg(\varphi_{d,\vec{x}}^*\mathcal{O}(1)|_{F_{n_1,\ldots,n_4}}) = \begin{cases}\min\{\dist(x_{n_i},\mathbb{Z})\} \text{ if } \sum_{i=1}^4\lfloor x_{n_i}\rfloor = d-1, \cr 0 \text{ otherwise, }\end{cases}\end{equation}
where $\dist(y,\mathbb{Z}) := \min\{y -\lfloor y \rfloor, \lceil y\rceil - y\}$.
\end{theorem}

\subsection{Proof of Theorem \ref{thm:ASgen}}
We proceed in several steps.  First, we verify the second line of (\ref{eq:ASgen}).  Next, we discuss the diagonal morphisms and Gale transform that allow for an inductive reduction to the case $n=4,d=1$.  We conclude by computing the polarization on $V_{1,4}\quotient_{\vec{y}}\SL_2\cong\mathbb{P}^1$ as a function of $\vec{y}\in\Delta(2,4)$.

\subsubsection{Degree zero}

\begin{lemma}
If $x_{n_i} \ge \alpha_i$, $i=1,\ldots,4$, for $\alpha_i \in\mathbb{Z}_{\ge 0}$ satisfying $\sum_{i=1}^4\alpha_i=d$, or if $x_{n_i} \le \beta_i$ for $\beta_i \in\mathbb{Z}_{\ge 1}$ with $\sum_{i=1}^4\beta_i=d+2$, then $F_{n_1,\ldots,n_4}$ is contracted by $\varphi_{d,\vec{x}}$.
\end{lemma}

\begin{proof}
This follows from the proofs of \cite[Propositions 4.1 and 4.2]{Giansiracusa}.
\end{proof}

Now $\sum_{i=1}^4 x_{n_i} = \sum_{j=1}^n x_j = d+1$, so \[\sum_{i=1}^4\lfloor x_{n_i}\rfloor \in \{d-2,d-1,d,d+1\}.\]  If this sum is $d$, then setting \[(\alpha_1,\ldots,\alpha_4)=(\lfloor x_{n_1} \rfloor, \ldots, \lfloor x_{n_4} \rfloor)\] and applying the preceding lemma shows that $\varphi_{d,\vec{x}}^*\mathcal{O}(1)$ has degree zero on $F_{n_1,\ldots,n_4}$.  If $\sum_{i=1}^4\lfloor x_{n_i}\rfloor = d+1$, then set $\alpha_i = \lfloor x_{n_i}\rfloor - 1$ for some $i$ and the same result applies.  On the other hand, if $\sum_{i=1}^4\lfloor x_{n_i}\rfloor = d-2$ then $\sum_{i=1}^4\lceil x_{n_i}\rceil = d+2$, so we can set \[(\beta_1,\ldots,\beta_4)=(\lceil x_{n_1} \rceil, \ldots, \lceil x_{n_4} \rceil)\] and apply the same lemma.  This verifies the second line of (\ref{eq:ASgen}).

\subsubsection{Diagonal morphisms}
Given a finite set $S=\{1,\ldots,k\}$ and a partition $S=\sqcup_{i=1}^l s_i$ with $|s_i| \ge 1$, there is a corresponding diagonal morphism $(\mathbb{P}^d)^l \rightarrow (\mathbb{P}^d)^k$.  This restricts to a morphism $V_{d,l} \rightarrow V_{d,k}$.  For $d\le k -3$ and $\vec{x}\in\Delta(d+1,k)$, let $x_{s_i} := \sum_{j\in s_i} x_j$.  If $x_{s_i} \le 1$ for all $i$, then $(x_{s_1},\ldots,x_{s_l})\in\Delta(d+1,l)$, and there is an induced morphism \[V_{d,l}\quotient_{(x_{s_1},\ldots,x_{s_l})}\SL_{d+1} \rightarrow V_{d,k}\quotient_{(x_1,\ldots,x_k)}\SL_{d+1}.\] This follows from the fact that GIT stability for these loci is determined by the amount of weight lying in various subspaces (cf. \cite[\S4.2]{Giansiracusa}).  This morphism pulls back the natural GIT polarization of the codomain to the GIT polarization of the domain.  If we omit the linearizations for such a map, then they are tacitly assumed to be of the form indicated above.

\begin{lemma}\label{lem:diag}
If $x_{s_i} \le 1$ for all $i=1,\ldots,l$, then there is a commutative diagram
\[\xymatrix{ \ovop{M}_{0,l} \ar[r]\ar[d] & \ovop{M}_{0,k} \ar[d] \\ V_{d,l}\quotient\SL_{d+1} \ar[r] & V_{d,k}\quotient\SL_{d+1}}\] where the top arrow sends a stable $l$-pointed curve to the stable $k$-pointed curve obtained by attaching maximally degenerate chains of $\mathbb{P}^1$s to the marked points according to the partition $\{1,\ldots,k\}=\sqcup_{i=1}^l s_i$.
\end{lemma}

\begin{proof}
It follows from the pointwise description of the vertical arrows \cite[\S4.1]{Giansiracusa}, the description of GIT stability \cite[\S4.2]{Giansiracusa}, and the condition $x_{s_i} \le 1$ that the chain sprouting out of the $i^{\text{th}}$ point gets contracted by the map $\ovop{M}_{0,k} \rightarrow V_{d,k}\quotient\SL_{d+1}$.
\end{proof}

\subsubsection{Gale transform}
The Gale transform is an involutive isomorphism \[(\mathbb{P}^d)^k\quotient_{\vec{x}}\SL_{d+1}~\widetilde{\rightarrow}~(\mathbb{P}^{k-d-2})^k\quotient_{\gamma(\vec{x})}\SL_{k-d-1}\] where, with the normalization assumption $\vec{x}=(x_1,\ldots,x_k)\in\Delta(d+1,k)$, we define \begin{equation}\label{eq:GaleWt}\gamma(\vec{x}) := (1-x_1,\ldots,1-x_k)\in\Delta(k-d-1,k).\end{equation}  This transformation was studied classically by Cobble when the points are distinct (so that no GIT quotients are necessary).  For GIT quotients with symmetric linearization this isomorphism was first proven in \cite[Corollary III.1]{DO88}.  For arbitrary weights, there are proofs in \cite{Hu05,Alp10}.

In \cite[\S6.2]{Giansiracusa} it is shown that the Gale transform induces an isomorphism \[\Gamma : V_{d,k}\quotient_{\vec{x}}\SL_{d+1}~\widetilde{\rightarrow}~V_{k-d-2,k}\quotient_{\gamma(x)}\SL_{k-d-1}.\]  If we omit the linearizations for this map, then they are assumed to be related by the involutive function $\gamma$ defined in (\ref{eq:GaleWt}).  As with the diagonal morphism, this Gale morphism preserves the natural GIT polarizations.

\begin{lemma}\label{lem:Gale}
For any $d \le k-3$ and $\vec{x}\in\Delta(d+1,k)$, the diagram
\[\xymatrix{\ovop{M}_{0,k} \ar[dr]^{~\varphi_{k-d-2,\gamma(\vec{x})}}\ar[d]_{\varphi_{d,\vec{x}}} & \\ V_{d,k}\quotient_{\vec{x}}\SL_{d+1} \ar[r]_{\Gamma~\text{ }~\text{ }~} & V_{k-d-2,k}\quotient_{\gamma(\vec{x})}\SL_{k-d-1}}\]
is commutative.
\end{lemma}

\begin{proof}
Since everything is separated, it is enough to check this on a dense open subset, so consider the following diagram: 
\[\xymatrix{\MM_{0,k} \ar[d]\ar[dr] & \\ U_{d,k}/\SL_{d+1} \ar[r] & U_{k-d-2,k}/\SL_{k-d-1}}\] Now $\MM_{0,k} \cong U_{1,k}/\SL_2$ and both arrows from it are defined by applying a Veronese map, so commutativity follows from \cite[Proposition III.2.3]{DO88}.
\end{proof}

\subsubsection{The main reduction}

\begin{proposition}\label{prop:reduction}
For any $d\le n-3$, $\vec{x}\in\Delta(d+1,n)$, and partition $\{1,\ldots,n\} = \sqcup_{i=1}^4 n_i$ such that $\sum_{i=1}^4\lfloor x_{n_i} \rfloor = d-1$, there is a commutative diagram
\[\xymatrix{ \ovop{M}_{0,4} \ar[r]\ar[d]_{\varphi_{1,\vec{y}}}^{\sim} & \ovop{M}_{0,n} \ar[d]^{\varphi_{d,\vec{x}}} \\ V_{1,4}\quotient_{\vec{y}}\SL_{2} \ar[r] & V_{d,n}\quotient_{\vec{x}}\SL_{d+1}}\]
where $\vec{y} = (y_1,y_2,y_3,y_4)$ with $y_i\in\{1,x_{n_i} - \lfloor x_{n_i} \rfloor, \lceil x_{n_i} \rceil - x_{n_i}\}$, such that the bottom arrow preserves GIT polarizations and the top arrow sends the fundamental class $[\ovop{M}_{0,4}]$ to $F_{n_1,\ldots,n_4}$.  Consequently, $\deg(\varphi_{d,\vec{x}}^*\mathcal{O}(1)|_{F_{n_1,\ldots,n_4}}) = \deg \mathcal{O}_{1,\vec{y}}(1)$.
\end{proposition}

This diagram is inspired by the one used in \cite{AS08} for the case $d=1$.

\begin{proof}
We first address the case $d=1$.  This is immediate from Lemma \ref{lem:diag}, since the condition $\sum_{i=1}^4 \lfloor x_{n_i} \rfloor = d - 1 = 0$ implies that $\lfloor x_{n_i} \rfloor = 0$ for all $i$, so $x_{n_i} \le 1$ and, moreover, $\vec{y} = (x_{n_1},\ldots,x_{n_4})$ is of the required form.  Next, consider the case $d=n-3$.  By Lemmas \ref{lem:diag} and \ref{lem:Gale}, there is a commutative diagram
\[\xymatrix{ \ovop{M}_{0,4} \ar[r]\ar[d] & \ovop{M}_{0,n} \ar[d] \ar[dr] & \\ V_{1,4}\quotient_{\vec{y}}\SL_{2} \ar[r] & V_{1,n}\quotient_{\gamma(\vec{x})}\SL_{2} \ar[r]_{\Gamma~} & V_{n-3,n}\quotient_{\vec{x}}\SL_{n-2}}\] with $\vec{y}=(|n_1| - x_{n_1},\ldots,|n_4| - x_{n_4})$, as long as $|n_i| - x_{n_i} \le 1$.  If $x_{n_i}\notin\mathbb{Z}$ for each $i$, then we claim that $|n_i| = \lceil x_{n_i} \rceil$.  Indeed, $x_j \le 1$ for each $j=1,\ldots, n$, so $x_{n_i} \le |n_i|$ and hence $\lceil x_{n_i} \rceil \le |n_i|$.  By assumption, $\sum_{i=1} \lfloor x_{n_i} \rfloor = d - 1 = n - 4$, hence \[\sum_{i=1}^4 \lceil x_{n_i} \rceil = \sum_{i=1}^4 (\lfloor x_{n_i} \rfloor + 1) = n = \sum_{i=1}^4 |n_i|,\] so this inequality is an equality.  If $x_{n_i}\in\mathbb{Z}$ for some $i$, then we instead deduce that $|n_i| = x_{n_i} + 1$.  In either case, the diagram exists and $\vec{y}$ is of the required form.

Having established these extreme cases, we now turn to the general case, for which we use induction on the number of marked points.  The base case $n=4$ is already settled, since it must have $d=1$.  So fix $N \ge 4$ and suppose the statement holds for all $n \le N$ and $d \le n-3$.  We must show that it also holds for $n=N+1$.  By the preceding paragraph, we can assume that $2 \le d \le n-4$.

If necessary, relabel the points so that the given partition $\{1,\ldots,n\} = \sqcup_{i=1}^4 n_i$ has $x_{n-1}$ and $x_n$ in the same part.  There are two cases to consider.  If $x_{n-1}+x_n \le 1$, then by setting $\vec{x}' := (x_1,\ldots,x_{n-2},x_{n-1}+x_n)$ we have a commutative diagram
\[\xymatrix{\ovop{M}_{0,n-1} \ar[r]\ar[d] & \ovop{M}_{0,n} \ar[d] \\ V_{d,n-1}\quotient_{\vec{x}'}\SL_{d+1} \ar[r] & V_{d,n}\quotient_{\vec{x}}\SL_{d+1}}\] by Lemma \ref{lem:diag}.  The result then follows from the inductive hypothesis, where we use the obvious partition of $\{1,\ldots,n-1\}$ derived from the original partition.  On the other hand, if $x_{n-1}+x_n > 1$, then  instead we consider the following diagram:
\[\xymatrix{ & \ovop{M}_{0,n-1} \ar[r]\ar[dl]\ar[d] & \ovop{M}_{0,n} \ar[dr] \ar[d] &   \\ V_{d-1,n-1}\quotient_{\vec{x}'}\SL \ar[r]_{\Gamma~} & V_{n-d-2,n-1}\quotient\SL \ar[r] & V_{n-d-2,n}\quotient\SL \ar[r]_{~\Gamma} & V_{d,n}\quotient_{\vec{x}}\SL}\] Here $\vec{x}' := (x_1,\ldots,x_{n-2},x_{n-1}+x_n -1)$.  The point is that $\gamma(\vec{x}) = (1-x_1,\ldots,1-x_n)$ and $(1-x_{n-1})+(1-x_n) = 2 - (x_{n-1}+x_n) < 1$, so the diagonal morphism in the bottom row exists.  Commutativity of this diagram then follows immediately from Lemmas \ref{lem:diag} and \ref{lem:Gale}, so as before we can apply the inductive hypothesis.
\end{proof}

\subsubsection{Polarizations on $\mathbb{P}^1$}
By Proposition \ref{prop:reduction}, it only remains to prove:

\begin{lemma}\label{lemma:AS}
For $\vec{y}=(y_1,\ldots,y_4)\in\Delta(2,4)$, the natural GIT polarization $\mathcal{O}(1)$ on $(\mathbb{P}^1)^4\quotient_{\vec{y}}\SL_2\cong\mathbb{P}^1$ has degree $\min\{y_1,\ldots,y_4,1-y_1,\ldots,1-y_4\}$.
\end{lemma}

\begin{remark}
This was first proven in \cite[Lemma 2.2]{AS08}.  We include an independent proof to illustrate that it can be achieved through a direct, elementary computation.  The proof in \cite{AS08} is more slick: it cleverly uses toric geometry to reduce the computation to that of determining the area of a triangle!
\end{remark}

\begin{proof}
Recall that any projective GIT quotient can be constructed as the projective spectrum of an invariant ring: \[X\quotient_L G = \Proj(\oplus_{r \ge 0}\HH^0(X, L^{\otimes r})^G).\]  In our case we set $L:=\mathcal{O}(y_1,\ldots,y_4)$, and then \[(\mathbb{P}^1)^4\quotient_{\vec{y}}\SL_2 = \Proj(\oplus_{r\ge 0}\HH^0((\mathbb{P}^1)^4,L)^{\SL_2}).\]  Therefore, as long as there are invariant global sections, we have \[\HH^0((\mathbb{P}^1)^4\quotient_{\vec{y}}\SL_2,\mathcal{O}(1)) = \HH^0((\mathbb{P}^1)^4,L)^{\text{SL}_2}.\] Since this quotient is isomorphic to $\mathbb{P}^1$, we conclude that \begin{equation}\label{eq:countdeg}\deg\mathcal{O}(1) = h^0(\mathcal{O}(1)) - 1 = \dim(\HH^0(L)^{\text{SL}_2}) - 1.\end{equation} Thus we are reduced to a straightforward problem in classical invariant theory: count the number of independent $\SL_2$-invariant global sections of $L$.  

If the $y_i$ are multiplied by a common factor then the degree of the corresponding polarization scales by that same factor, so to compute this degree we can assume that $y_i\in\mathbb{N}$ and $y := \sum_{i=1}^4 y_i\in 2\mathbb{N}$.  This guarantees the existence of invariant sections.  If we think of coordinates on $(\mathbb{P}^1)^4$ as a $2\times 4$ matrix, then $\SL_2$-invariant polynomials are given by the $2\times 2$ minors of this matrix.  We denote the minor with columns $i$ and $j$ by the $2\times 1$ tableaux \[\left.\begin{array}{|c|}\hline i \\\hline j \\\hline \end{array}\right.\]  Horizontal juxtaposition of tableaux indicates the product of the corresponding minors.  It is a classical fact that a basis for the invariant functions we are interested in is then given by all $2\times \frac{y}{2}$ tableaux with entries in $\{1,2,3,4\}$ that are semi-standard (i.e., the entries increase strictly down the columns and weakly across the rows) such that the number $i$ appears exactly $y_i$ times.  To count how many such diagrams are possible, we can label the weights so that $y_1 \ge y_2 \ge y_3 \ge y_4$.  Moreover, if $y_1 > \frac{y}{2}$, then the semistable locus is empty, so we can assume $y_1 \le \frac{y}{2}$.

First consider the case that $y_1+y_4 = \frac{y}{2}$.  Here the possible diagrams are 
\[\left.\begin{array}{|c|c|c|c|c|c|c|c|c|}\hline 1 & \cdots & 1 & 1 & \cdots & 2 & \cdots & 2 & 2\\\hline \cdots & 2 & 3 & 3 & \cdots & 4 & \cdots & 4 & 4\\\hline \end{array}\right.,\]

\[\left.\begin{array}{|c|c|c|c|c|c|c|c|c|}\hline 1 & \cdots & 1 & 1 & \cdots & 2 & \cdots & 2 & \mathbf{3} \\\hline \cdots & 2 & \mathbf{2} & 3 & \cdots & 4 & \cdots & 4 & 4 \\\hline \end{array}\right.,\] 
\[\vdots\]
\[\left.\begin{array}{|c|c|c|c|c|c|}\hline 1 & \cdots & 1 & 3 & \cdots & 3 \\\hline 2 & \cdots & 3 & 4 & \cdots & 4 \\\hline \end{array}\right.\]
That is to say, one can start by filling in all the entries in order, then the only other tableaux are obtained by taking a 3 from the second row and switching it with a 2 in the first row.  We obtain exactly $y_4 + 1 = \frac{y}{2} - y_1 + 1$ diagrams in this fashion.  If $y_1+y_4 > \frac{y}{2}$, then the diagrams look similar, except now the 1s in the first row can lie above the 4s of the second row.  The number of 2s in the first row that we can swap with 3s in the second row is $\frac{y}{2} - y_1$, so there are $\frac{y}{2}-y_1+1$ total diagrams possible.  Thus we have $\frac{y}{2}-y_1 +1$ diagrams whenever $y_1+y_4\ge\frac{y}{2}$.

For $y_1+y_4 < \frac{y}{2}$ a similar argument applies, the only difference is that the 1s in the top row end before the 4s in the bottom row begin.  Because we cannot have a 3 on the top row above a 3 in the bottom row, the only place we can put 3s in the top row is above the 4s, so the number of possible diagrams is $y_4+1$.  

Scaling back down to fractional weights via the original normalization $y=2$, we see from (\ref{eq:countdeg}) that the degree is $1-y_1$ if $1 - y_1 \le y_4$, and it is $y_4$ otherwise.  Since the weights were ordered $y_1 \ge \cdots \ge y_4$, this is equivalent to the claimed formula.
\end{proof}

This completes the proof of Theorem \ref{thm:ASgen}, which is our main technical tool.  

\section{$\mathcal{CB}(A,1)$ and GIT}\label{CBGIT}
In this section we prove Theorem \ref{thm:morphism}, which says that every nontrivial type $A$, level one, conformal blocks divisor induces a morphism $\ovop{M}_{0,n} \rightarrow V_{d,k}\quotient_{\vec{x}}\SL_{d+1}$.

\subsection{Setup}
For $\mathfrak{sl}_m$, the vector of conformal blocks weights is of the form $w=(\omega_{c_1},\ldots, \omega_{c_n})$, where the $\omega_{c_i}$ are fundamental dominant weights with $0 \le c_i \le m$.  Let $\vec{c} := (c_1,\ldots,c_n)$, and denote the corresponding divisor by $\mathbb{D}^{\mathfrak{sl}_m}_{1,\vec{c}} := \mathbb{D}^{\sL_m}_{1,(\omega_{c_1},\ldots,\omega_{c_n})}$.

By ``propagation of vacua'' (cf. \cite[Proposition 2.4.(1)]{Fak}), if $c_i$ is 0 or $m$ for some $i$, then the induced morphism factors through the map $\ovop{M}_{0,n} \rightarrow \ovop{M}_{0,n-1}$ given by forgetting the $i^{\text{th}}$ marked point, so assume for the remainder of this section that $1 \le c_i \le m-1$ for each $i$.  

The following lemma and proof were kindly communicated to us by Fakhruddin.

\begin{lemma}\label{lem:Fakh}
For $\mathbb{D}^{\mathfrak{sl}_m}_{1,\vec{c}}$ to be nontrivial, it must satisfy $\sum_{i=1}^n c_i = m(d+1)$ for some $d\in\{1,\ldots,n-3\}$.
\end{lemma}

\begin{proof}
By the theory of conformal blocks, in order to get a nonzero bundle the sum of the weights must lie in the root lattice, so $m | (\sum_{i=1}^n c_i)$, hence $\sum_{i=1}^n c_i=m(d+1)$ for some $d \ge 0$.  It follows immediately from \cite[Proposition 5.2]{Fak} that if $d=0$ then the bundle is trivial, so we have $d \ge 1$.  On the other hand, there is a duality \[\mathbb{D}^{\mathfrak{sl}_m}_{1,(c_1,\ldots,c_n)} = \mathbb{D}^{\mathfrak{sl}_m}_{1,(m-c_1,\ldots,m-c_n)}\] arising from the representation theory of $\mathfrak{sl}_m$, so we must have $\sum_{i=1}^n (m-c_i) = m(d'+1)$ for some $d' \ge 1$.  Thus \[m(d+1) = \sum_{i=1}^n c_i = m(n-d'-1),\] and hence $d = n - d' - 2 \le n-3$.
\end{proof}

\subsection{Relation to GIT}\label{section:RelGIT}
The first correspondence between conformal blocks and GIT quotients was found by Fakhruddin, where he showed in \cite[Theorem 4.5]{Fak} that Kapranov's morphisms $\ovop{M}_{0,n} \rightarrow (\mathbb{P}^1)\quotient_{\vec{x}}\SL_2$ are all induced by $\mathfrak{sl}_2$ conformal blocks divisors.  By \cite[Remark 5.3]{Fak}, these divisors also occur as elements of $\mathcal{CB}(A,1)$.  It is shown in \cite[Theorem 1.2]{Giansiracusa} that the morphisms $\ovop{M}_{0,n} \rightarrow V_{d,n}\quotient_{\vec{x}}\SL_{d+1}$ for $S_n$-invariant linearization $\vec{x}$ are induced by elements of $\mathcal{CB}(A,1)$ corresponding to $\mathfrak{sl}_n$ bundles.  In personal correspondences, Fakhruddin suggested Lemma \ref{lem:Fakh} and observed that the integers $d$ lie in the same range as the $d$ in the GIT quotients $V_{d,n}\quotient\SL_{d+1}$, hinting at the following:  

\begin{theorem}\label{thm:CBGIT}
If $\mathbb{D}^{\mathfrak{sl}_m}_{1,\vec{c}}$ is nontrivial, so that $\sum_{i=1}^n c_i = m(d+1)$ for $1\le d \le n-3$, then its divisor class coincides with that of $\varphi_{d,\vec{c}}^*\mathcal{O}(1)$.
\end{theorem}

\begin{proof}
It is enough to show that both classes have the same intersection number with every F-curve.  For the GIT polarizations this formula is given by Theorem \ref{thm:ASgen}, and for the conformal blocks bundles it is given by Fakhruddin's formula \cite[Proposition 5.2]{Fak}.  The formula for GIT quotients is stated for normalized linearizations that lie in the hypersimplex, so it is convenient to write $\vec{x} := (\frac{c_1}{m},\ldots,\frac{c_n}{m})$, since then $\vec{x}\in\Delta(d+1,n)$.  Since $\mathcal{O}_{d,\vec{c}}(1) = \mathcal{O}_{d,\vec{x}}(m)$, we must show that for any F-curve class $\{1,\ldots,n\}=\sqcup_{i=1}^4 n_i$ there is an equality \[\deg(\mathbb{D}^{\mathfrak{sl}_m}_{1,\vec{c}}|_{F_{n_1,\ldots,n_4}}) = m\deg(\varphi_{d,\vec{x}}^*\mathcal{O}(1)|_{F_{n_1,\ldots,n_4}}).\] 

If we write \[\nu_i := \sum_{i\in n_i} c_i~(\text{mod}~m) \in \{0,\ldots, m-1\},\] then part of Fakhruddin's formula is that $\sum_{i=1}^4 \nu_i \ne 2m \Rightarrow \deg(\mathbb{D}^{\mathfrak{sl}_m}_{1,\vec{c}}|_{F_{n_1,\ldots,n_4}}) = 0$.  We claim that \begin{equation}\label{eq:equiv}\sum_{i=1}^4 \nu_i = 2m \Longleftrightarrow \sum_{i=1}^4 \lfloor x_{n_i} \rfloor = d - 1.\end{equation} Indeed, $\sum_{i=1}^4 x_{n_i} = d+1$, so $\sum_{i=1}^4 \lfloor x_{n_i} \rfloor = d - 1 \Leftrightarrow \sum_{i=1}^4 (x_{n_i} - \lfloor x_{n_i} \rfloor) = 2$.  On the other hand, $\nu_i = mx_{n_i}~(\text{mod}~m)$, so $\frac{\nu_i}{m} = x_{n_i}~(\text{mod}~1) = x_{n_i} - \lfloor x_{n_i} \rfloor$.  Thus $\sum_{i=1}^4\nu_i = 2m \Leftrightarrow 2 = \sum_{i=1}^4 \frac{\nu_i}{m} = \sum_{i=1}^4 (x_{n_i} - \lfloor x_{n_i} \rfloor)$.  This verifies the claim.

For the remainder of the proof, assume that the equivalent conditions in (\ref{eq:equiv}) are satisfied. In this case, writing $\nu_{\max} = \max\{\nu_1,\ldots,\nu_4\}$ and $\nu_{\min} = \min\{\nu_1,\ldots,\nu_4\}$, Fakhruddin's formula says \[\deg(\mathbb{D}^{\mathfrak{sl}_m}_{1,\vec{c}}|_{F_{n_1,\ldots,n_4}}) = \min\{\nu_{\min},m-\nu_{\max}\},\] Now $0\le \nu_i < m$ and $\frac{\nu_i}{m} = x_{n_i} - \lfloor x_{n_i} \rfloor$, so 
\begin{eqnarray*}
\min\{\nu_{\min},m-\nu_{\max}\}&=&\min\{\nu_1,\ldots,\nu_4,m-\nu_1,\ldots, m-\nu_4\}\\
& = & m\cdot \min\{\frac{\nu_1}{m},\ldots,1-\frac{\nu_1}{m},\ldots\}\\
& = & m\cdot \min\{x_{n_1} - \lfloor x_{n_1} \rfloor, \ldots,\lceil x_{n_1} \rceil - x_{n_1}, \ldots\}\\
& = & m\cdot\min\{\dist(x_{n_i},\mathbb{Z})\}.
\end{eqnarray*}
This completes the proof.
\end{proof}

\section{Finite Generation of $\mathcal{CB}(A,1)$}\label{finitegeneration}
In this section we prove Theorem \ref{thm:finite}, which says that $\mathcal{CB}(A,1)$ spans a finitely generated cone.  A first step is to prove a converse to Theorem \ref{thm:morphism}.

\begin{proposition}\label{prop:converse}
Any birational morphism $\ovop{M}_{0,n} \rightarrow V_{d,n}\quotient_{\vec{x}}\SL_{d+1}$ extending the obvious identification of interiors is induced by an element of $\mathcal{CB}(A,1)$.
\end{proposition}

\begin{proof}
By separatedness, any such morphism is given by the morphism $\varphi_{d,\vec{x}}$, so we just need to find a conformal blocks divisor inducing this morphism.  We can assume that $\vec{x}\in\Delta(d+1,n)$.  Then there is $m\in\mathbb{Z}$ that clears the denominators: setting $c_i := mx_i$ yields $\vec{c}=(c_1,\ldots,c_n)\in\mathbb{Z}^n$.  The hypothesis that $\vec{x}$ lies in the hypersimplex then implies that $\sum_{i=1}^n c_i = m(d+1)$.  Now $\varphi_{d,\vec{x}}^*\mathcal{O}(1) = \frac{1}{m}\cdot \varphi_{d,\vec{c}}^*\mathcal{O}(1)$, and by Theorem \ref{thm:CBGIT} the class of $\varphi_{d,\vec{c}}^*\mathcal{O}(1)$ coincides with the class of $\mathbb{D}^{\mathfrak{sl}_m}_{1,\vec{c}}$.
\end{proof}

If the linearization $\vec{x}$ contains any zeros, then the morphism from $\ovop{M}_{0,n}$ to the GIT quotient is not birational, but it factors as a forgetful map $\ovop{M}_{0,n} \rightarrow \ovop{M}_{0,k}$ followed by a birational morphism to the GIT quotient.  Thus, by combining Theorem \ref{thm:morphism} and Proposition \ref{prop:converse}, we obtain the following:

\begin{corollary}\label{cor:cones}
The cone spanned by $\mathcal{CB}(A,1)$ coincides with the cone spanned by \[\{\varphi_{d,\vec{x}}^*\mathcal{O}(1) : d\in\{1,\ldots,n-3\},\vec{x}\in\Delta(d+1,n)\}.\]
\end{corollary}

We will use this description to prove finite generation.

\subsection{GIT cones}

Alexeev and Swinarski define in \cite{AS08} a subcone of $\Nef(\M_{0,n})$, called the GIT cone, as follows.  For each linearization $\vec{x}\in\Delta(2,n)$, one has the morphism $\varphi_{1,\vec{x}} : \M_{0,n} \rightarrow (\mathbb{P}^1)^n\quotient_{\vec{x}}\SL_2$, first introduced by Kapranov \cite{KapChow}.  One obtains a nef divisor $\varphi_{1,\vec{x}}^*\mathcal{O}(1)$ on $\M_{0,n}$, and the convex cone spanned by these divisors is the GIT cone.  Since $(\mathbb{P}^1)^n=V_{1,n}$, we shall call this the \emph{degree one} GIT cone.  By fixing any $d\in\{1,\ldots,n-3\}$, we can define an analogous cone, the degree $d$ GIT cone, by pulling back the GIT polarization on the quotients $V_{d,n}\quotient_{\vec{x}}\SL_{d+1}$ for all $\vec{x}\in\Delta(d+1,n)$.  Thus Corollary \ref{cor:cones} says that $\mathcal{CB}(A,1)$ is spanned by all the GIT cones.  In fact, by Gale duality, the degree $d$ GIT cone coincides with the degree $n-d-2$ GIT cone (cf. \cite[\S6]{Giansiracusa}), so it is enough to have $1 \le d \le \lfloor \frac{n}{2} \rfloor - 1$.

To conclude the proof of Theorem \ref{thm:finite}, it only remains to prove the following:

\begin{proposition}
For each $d$, the degree $d$ GIT cone is finitely generated.
\end{proposition}

\begin{proof}
By \cite[Theorem 2.3]{Tha96}, GIT polarizations vary linearly within each GIT chamber, so the function $\theta : \Delta(d+1,n) \rightarrow \N^1(\ovop{M}_{0,n})$ defined by $\theta(\vec{x}) := \varphi^*_{d,\vec{x}}\mathcal{O}(1)$ is piecewise linear.  Moreover, by the same theorem, as the linearization moves from the interior of a chamber to an adjacent wall, there is an induced morphism of quotients. So $\theta$ is also continuous.  This implies that any extremal ray for the cone generated by the image of $\theta$ must come from a vertex of the GIT chamber decomposition of $\Delta(d+1,n)$.  By \cite[Theorem 2.4]{Tha96}, there are only finitely many walls, hence finitely many vertices, hence finitely many extremal rays.
\end{proof}

\begin{remark}
Not only does the above proof show that $\mathcal{CB}(A,1)$ has finitely many extremal rays, but at least in theory one can find them all.  Indeed, each extremal ray for this cone comes from an extremal ray for one of the GIT cones, so it is enough to find the extremal rays of each GIT cone.  For a fixed value of $d$, the chamber decomposition of $\Delta(d+1,n)$ is cut out by all walls of the form $\sum_{i\in I}x_i = k$ for $I\subseteq\{1,\ldots,n\}$ and $1\le k \le d$ \cite[Example 3.3.24]{DH98}.  So if one could solve the problem in polyhedral geometry of determining all vertices for this decomposition, then one could use Theorem \ref{thm:ASgen} to compute the corresponding divisor classes and obtain all possible extremal rays in this fashion.  Unfortunately, determining the vertices from the description of the walls seems, however elementary, to be quite difficult.  See \cite{AS08}, for example, which studies the case $d=1$. 
\end{remark}

\begin{remark}\label{D}
Fakhruddin has shown that the closure of the cone generated by $\mathcal{CB}(D,1)$ is finitely generated \cite[5.2.7, Prop 5.6]{Fak}, so such divisors yield at most finitely many extremal rays of the nef cone.  The exceptional groups provide only finitely many conformal blocks bundles of a given level.  In fact, Fakhruddin proved that $\mathcal{CB}(\op{F_4},1)$ and $\mathcal{CB}(\op{G_2},1)$ are contained in the ample cone \cite[5.2.8]{Fak}, so there are no extremal rays of $\op{Nef}(\ovop{M}_{0,n})$ given by such divisors.  Moreover, $\mathcal{CB}(\mathfrak{e}_6,1)\cup\mathcal{CB}(\mathfrak{e}_7,1)\subseteq\mathcal{CB}(A,1)$ and $\mathcal{CB}(\mathfrak{e}_8,1)=\{0\}$ \cite[5.2.5-7]{Fak}.  The only level one conformal blocks divisor not known to yield only finitely many extremal rays of the nef cone are of type $B$ and $C$.
\end{remark}

\section{Identitities and applications}\label{section:Identities}
The cone generated by $\mathcal{CB}(A,1)$ has only finitely many extremal rays, by Theorem \ref{thm:finite}, yet there are infinitely many type $A$ Lie algebras and choices of weights.  So one expects that many different Lie algebras and weights will yield proportional conformal blocks divisor classes.  In this section, in Proposition \ref{identities}, we explicitly describe an instance of this.   As an application, we consider $S_n$-invariant elements of $\mathcal{CB}(A,1)$.  In \cite[Theorem B]{agss} a family of $\lfloor \frac{n}{2} \rfloor-1$ elements of $\mathcal{CB}(\sL_n,1)$ were shown to be extremal in the nef cone of $\ovop{M}_{0,n}/S_n$, and in Corollary \ref{extremal} we show that these are the only extremal rays of the symmetric nef cone that come from $\mathcal{CB}(A,1)$.   Two of these extremal rays induce maps from the hyperelliptic and cyclic trigonal loci in $\ovop{M}_g$ to Satake's compactification of $\op{A}_g$, and in Corollary \ref{SatakeImage} we use Theorem \ref{thm:morphism} to describe their images as GIT quotients.

\subsection{The identities}
We continue to use the notation $\mathbb{D}^{\mathfrak{sl}_m}_{1,\vec{c}} := \mathbb{D}^{\sL_m}_{1,(\omega_{c_1},\ldots,\omega_{c_n})}$, where $\vec{c}=(c_1,\ldots,c_n)$.   

\begin{proposition}\label{identities}
For any  integer $k\ge 1$, there is an equality of divisor classes 
$$\mathbb{D}^{\sL_m}_{1,\vec{c}}=\frac{1}{k}\mathbb{D}^{\sL_{km}}_{1,k\vec{c}}.$$
\end{proposition}

\emph{Combinatorial Proof}.  It suffices to show that $\mathbb{D}^{\sL_m}_{1,\vec{c}}$ and $\frac{1}{k}\mathbb{D}^{\sL_{km}}_{1,k\vec{c}}$
have the same degree on each $\op{F}$-curve.   This follows immediately from \cite[Proposition $5.2$]{Fak}. \hfill $\Box$

\emph{Geometric Proof}. 
By Lemma \ref{lem:Fakh}, we must have $\sum_{i=1}^n c_i = (d+1)m$ for some $1\le d \le n-3$.  But then $\sum_{i=1}^n kc_i = (d+1)km$, so applying Theorem \ref{thm:CBGIT} twice yields \[\mathbb{D}^{\sL_m}_{1,\vec{c}}=\varphi_{d,\vec{c}}^*\mathcal{O}(1) = \frac{1}{k}\varphi_{d,k\vec{c}}^*\mathcal{O}(1) = \frac{1}{k}\mathbb{D}^{\sL_{km}}_{1,k\vec{c}},\] as desired. \hfill $\Box$

\begin{remark}
This proposition yields identities for conformal blocks divisors given by Lie algebras besides those of type $A$.  Indeed, Fakhruddin has shown that level one divisors for $\mathfrak{e}_6$ are the same as for $\mathfrak{sl}_3$, and level one divisors for $\mathfrak{e}_7$ are the same as for $\mathfrak{sl}_2$ \cite[5.2.5-6]{Fak}.
\end{remark}
  
\subsection{Extremal rays of the symmetric cone}

We say that a divisor $\mathbb{D}^{\sL_m}_{1, \vec{c}}$ is \emph{symmetric} if $\vec{c}=(j, j, \ldots, j)$ consists of an $S_n$-invariant set of weights. 
In this case, the class of the divisor is invariant under the action of $S_n$ on $\ovop{M}_{0,n}$ defined by permuting the labeling of the points.  One may regard the symmetric conformal blocks divisors as elements of the nef cone of $\ovop{M}_{0,n}/S_n$.  By \cite{GKM}, the birational geometry of $\ovop{M}_{0,n}/S_n$ is intimately connected to that of $\ovop{M}_n$.  For example, if one knew all the extremal rays of the nef cone of $\ovop{M}_{0,n}/S_n$, then one would know the extremal rays of $\Nef(\ovop{M}_n)$.  

In \cite{agss}, the authors showed that the divisors $\{\mathbb{D}^{\sL_n}_{1,\{j,\ldots,j\}}: 2 \le j \le \lfloor \frac{n}{2} \rfloor \}$ generate distinct extremal rays of the symmetric nef cone.   Using Proposition \ref{identities}, we can identify conformal blocks divisors given by infinitely many different choices of Lie algebras and sets of weights that generate extremal rays of $\op{Nef}(\ovop{M}_{0,n}/\op{S}_n)$:

\begin{corollary}\label{extremal}
Every nonzero divisor
$\mathbb{D}^{\sL_m}_{1, \{j,\ldots,j\}}$
generates an extremal ray of $\op{Nef}(\ovop{M}_{0,n}/\op{S}_n)$.     There are only $\lfloor \frac{n}{2} \rfloor -1$ rays of this form.
\end{corollary}

\begin{proof}By Lemma \ref{lem:Fakh}, the divisor
$\mathbb{D}^{\sL_m}_{1, \{j,\ldots,j\}}$ is nonzero if and only if $nj=m(d+1)$, for some $d+1 \in \{1,\ldots,n-2\}$.   By Proposition \ref{identities}:
\begin{multline}
\mathbb{D}^{\sL_m}_{1, \{j,\ldots, j\}}=\frac{1}{d+1}\mathbb{D}^{\sL_{m(d+1)}}_{1, \{j(d+1),\ldots, j(d+1)\}}\\
=\frac{1}{d+1}\mathbb{D}^{\sL_{nj}}_{1, \{j(d+1),\ldots, j(d+1)\}}=\frac{j}{d+1}\mathbb{D}^{\sL_{n}}_{1, \{d+1,\ldots, d+1\}}.
\end{multline}
By \cite[Theorem B]{agss}, the divisor $\mathbb{D}^{\sL_{n}}_{1, \{d+1,\ldots, d+1\}}$ is an extremal ray of the symmetric nef cone.
In particular, all such nonzero symmetric divisors are of the form $\mathbb{D}^{\sL_{n}}_{1, \{k,\ldots, k\}}$, for $1\le k \le n$. By Gale duality, $\mathbb{D}^{\sL_{n}}_{1, \{k,\ldots, k\}}=\mathbb{D}^{\sL_{n}}_{1, \{n-k,\ldots, n-k\}}$, and by \cite{ags}, $\mathbb{D}^{\sL_{n}}_{1, \{1,\ldots, 1\}}=0$, 
and the remaining $\lfloor \frac{n}{2}\rfloor  -1$ are distinct. 

\end{proof}

\subsection{The hyperelliptic and cyclic trigonal loci}
Theorem \ref{thm:CBGIT} can be used to study certain geometric loci in Satake's compactification of the moduli space of principally polarized abelian varieties.  As we show, the images of the hyperelliptic and cyclic trigonal loci under the extended Torelli map $\ovop{M}_g \rightarrow \ovop{A}_g^{Sat}$ can be constructed as GIT quotients parametrizing configurations of points lying on quasi-Veronese curves.  To motivate this, we first provide some background.

A hyperelliptic curve is a smooth projective curve that admits a $2:1$ map to  $\mathbb{P}^1$.  A cyclic trigonal curve is a smooth projective curve that can be realized as a degree $3$ cyclic cover of $\mathbb{P}^1$.   The set of hyperelliptic and cyclic trigonal curves form natural loci in $\op{M}_g$, and their closures in $\ovop{M}_g$ often reveal important information about $\ovop{M}_g$ itself.  For example, a crucial ingredient in Cornalba-Harris' celebrated result about the nef cone of $\ovop{M}_g$ \cite[Theorem 1.3]{CornalbaHarris} is a description of the Picard group of the hyperelliptic locus.   More recently, the intersection theory on the moduli space of stable hyperelliptic curves has been shown to be closely tied to the steps of the log minimal model program for $\ovop{M}_g$ \cite{HyeonLeeHyperelliptic}.

The Hodge class $\lambda$ induces the extended Torelli morphism $f_\lambda : \ovop{M}_g \rightarrow \ovop{A}_g^{Sat}$.  There is a well-known isomorphism $h$ from $\ovop{M}_{0,2(g+1)}/S_{2(g+1)}$  to the closure $\ovop{H}_g$ of the hyperelliptic locus in $\ovop{M}_g$, which associates to a $2(g+1)$-pointed
rational curve $C$, the double cover of $C$ branched at its (unordered) marked points \cite{AL}.  From \cite{CornalbaHarris}, one can write down the class of $h^*(\lambda)$, and therefore the divisor that gives the composition $h \circ f_{\lambda}$ (cf. \cite[Lemma 7.1]{ags}).   In  \cite{FedorchukCyclic},  Fedorchuk generalizes this map to cyclic $p$-covering morphisms
$f_{n,p} : \ovop{M}_{0,n} \rightarrow \ovop{M}_g$, for $g=\frac{(n-2)(p-1)}{2}$, by associating to each $n$-pointed rational curve, its cyclic degree $p$ cover branched over the marked points.  Moreover, he expresses the class of $f_{n,p}^*(\lambda)$ as a sum of conformal blocks divisors.   
\begin{definition}
Let  $$\ovop{H}_g = h(\ovop{M}_{0,2(g+1)}/S_{2(g+1)}) \subseteq \ovop{M}_g$$ be the stable hyperelliptic locus, and for $g$ such that $3|(g+2)$, let $$\ovop{CT}_g = f_{g+2,3}(\ovop{M}_{0,g+2}/S_{g+2})\subseteq \ovop{M}_g$$ be the stable cyclic trigonal locus.
\end{definition}

\begin{proposition}\label{SatakeImage} The normalization of $f_{\lambda}(\ovop{H}_g)$ is isomorphic to the normalization of the GIT quotient $(V_{g,2g+2}\quotient\SL_{g+1})/S_{2g+2}$ parameterizing configurations of $2g+2$ unordered points in $\mathbb{P}^{g}$ supported on a quasi-Veronese curve.  The normalization of $f_{\lambda}(\ovop{CT}_g)$ is isomorphic to the normalization of the quotient $(V_{\frac{g+2}{3},g+2}\quotient\SL_{\frac{g+5}{3}})/S_{g+2}$ parameterizing $g+2$ points in $\mathbb{P}^{\frac{g+2}{3}}$ on a quasi-Veronese curve.
\end{proposition}

\begin{proof} By \cite[Theorem 7.2]{ags} and Theorem \ref{thm:CBGIT} we have $$h^*(\lambda)= \frac{1}{2}\mathbb{D}^{\sL_{2}}_{1,(1,\ldots,1)} = \frac{1}{2}\varphi^*_{g,(1,\ldots,1)}\mathcal{O}(1),$$ from which the first statement follows.  For $g$ such that $3 | (g+2)$, by \cite[Theorem 4.4]{FedorchukCyclic}, one has 
$$f_{g+2,3}^*(\lambda)= \frac{1}{3} \mathbb{D}^{\sL_{g+2}}_{1,(\frac{g+2}{3},\ldots,\frac{g+2}{3})}+\frac{1}{3} \mathbb{D}^{\sL_{g+2}}_{1,(\frac{2(g+2)}{3},\ldots,\frac{2(g+2)}{3})},$$
 but by duality, $\mathbb{D}^{\sL_{g+2}}_{1,(\frac{2(g+2)}{3},\ldots,\frac{2(g+2)}{3})}=\mathbb{D}^{\sL_{g+2}}_{1,(\frac{g+2}{3},\ldots,\frac{g+2}{3})}$,  so $f^*_{g+2,3}(\lambda) = \frac{2}{3}\mathbb{D}^{\sL_{g+2}}_{1,(\frac{g+2}{3},\ldots,\frac{g+2}{3})}$, which by Theorem \ref{thm:CBGIT} is proportional to $\varphi^*_{\frac{g+2}{3},(1,\ldots,1)}\mathcal{O}(1)$.
\end{proof}

\end{document}